\documentclass[a4paper,11pt,draft]{article}
\usepackage[english]{babel}
\usepackage{amsmath,amsfonts,amsthm}
\usepackage{amssymb}
\usepackage[cp1251]{inputenc}
\newtheorem{theorem}{\sc\bf Theorem.}
\newtheorem{lemma}{\sc\bf Lemma}
\title{Compact spacelike surfaces in the 3-dimensional de Sitter space. }
\author{A.A. Borisenko \footnote{suported by research grant DFFD of Ukrainian Ministry of Education
and Science, No 01. 07/ 00132.}}
\date{}
\begin{document}
\maketitle{
 Kharkov National University, Faculty of Mathematics
and Mechanics,

Geometry Department, Svobodi sq., 4, Kharkov, 61077, Ukraine

e-mail: borisenk@univer. kharkov.ua}
 \vspace{0.5cm}

In this paper we establish several sufficient conditions for a
compact spacelike surface in the 3-dimensional de Sitter space to
be totally geodesic or spherical.

{\bf Key words:} De Sitter space, Compact spacelike surface,
Second fundamental form, Gaussian curvature; Totally umbilical
round sphere

{\bf Mathematics Subject Classification (2000) }  primary 53C42;
secondary 53B30, 53C45
 \vspace{0.5cm}

Let $E_1^4$  be a  $4$-dimensional Lorentz-Minkowski space, that
is, the space  endowed with the Lorentzian metric tensor
$\langle\, ,\rangle$ given by
$$
\langle ,\rangle=(dx_1)^2+(dx_2)^2+(dx_3)^2-(dx_0)^2,
$$
where $(x_1, x_2, x_3, x_0)$  are the canonical coordinates of
$E_1^4$. The $3$-dimensional unitary de Sitter space is defined as
the following hyperquadratic of $E_1^4$.
$$
S_1^3=\{x\in R^4: \langle x, x\rangle=1\}
$$
As it is well known, $S_1^3$ inherits from $E_1^4$ a
time-orientable Lorentzian metric which makes it the standard
model of a Lorentzian space of constant sectional curvature one. A
smooth immersion $\psi\colon F\to S_1^3\subset E_1^4$ of a
$2$-dimensional connected manifold $M$ is said to be a spacelike
surface if the induced metric via $\psi$ is a Riemannian metric on
$M$, which, as usual, is also denoted by $\langle \,, \rangle$.
The time-orientation of $S_1^3$ allows us to define a (global)
unique timelike unit normal field $n$ on $F$, tangent to $S_1^3$,
and hence we may assume that $F$ oriented by $n$. We will refer to
$n$ as the Gauss map of $F$.

We note  that Lobachevsky space $L^3$ is the set of points
$$
L^3=\{x\in E_1^4: \langle x, x\rangle=-1, x_0>0\}.
$$

 It is well known that  a compact spacelike surface in the
$3$-dimensional de Sitter space $S_1^3$ is diffeomorphic to a
sphere $S^2$. Thus, it is interesting to look for additional
assumptions for such a surface to be totally geodesic or totally
umbilical round sphere.

There are  two possible kinds of geometric assumptions: extrinsic,
that is re\-la\-ti\-ve to the second fundamental form, and
intrinsic, namely, concerning to the Gaussian curvature of the
induced metric. As regards to the extrinsic approach, Ramanathan
\cite{Ram} proved that every compact spacelike surface in $S_1^3$
of constant mean curvature is totally umbilical. This result was
generalized to hypersurface of any dimension by Montiel
\cite{Mon}. J.Aledo and A.Romero characterize the compact
spacelike surfaces in $S_1^3$ whose second fundamental form
defines a Riemannian metric. They studied the case of constant
Gaussian curvature $K_{II}$ of the second fundamental form,
proving that the totally umbilical round spheres are the only
compact spacelike surfaces in $S_1^3$ with $K<1$ and constant
$K_{II}$ \cite{AledRom}. With respect to the intrinsic approach Li
\cite{Li} obtained that compact spacelike surface of constant
Gaussian curvature is totally umbilical. And he proved there is no
complete spacelike surface in $S_1^3$ with constant Gaussian
curvature $K>1$. J.Aledo and A.Romero proved the same result
without condition that Gaussian curvature is constant
\cite{AledRom}. But it is true more general result.
\begin{theorem}\label{T1}
 Let $F$ be a $\mathcal{C}^2$-regular complete spacelike surface in de
Sitter space $S_1^3$. If Gaussian curvature $K\geqslant 1$ then
the surface $F$ is totally geodesic great sphere with Gaussian
curvature $K=1$.
\end{theorem}
S. N. Bershtein proved that an explicitly given saddle surface
over a whole plane in the Euclidean space $E^3$ with slower than
linear
 growth at infinity must be a cylinder. He proved this theorem for
 surfaces of class $\mathcal{C}^2$ \cite{Ber1}, and it was
 generalized to the non-regular case in \cite{AV}.

A surface $F^2$ of smoothness class $C^1$ in $S^3$ may be
projected univalently into a great sphere $S^2_0$ if the great
spheres tangent to $F^2$ do not pass through points $Q_1$, $Q_2$
polar to $S^2_0$.

The surface $F^2$ in $S^3$ is called a saddle surface if any
closed rectifiable contour $\cal L$, that is in the intersection
of $F^2$  with an arbitrary great sphere $S^2$ in $S^3$, lies in
an open hemisphere, and is deformable to a point in the surface
can be spanned by a two-dimensional simply connected surface $Q$
contained in $F^2\cap S^2$. In other words, from the surface it is
impossible to cut off a crust by a great sphere $S^2$, that is, on
$F^2$ there do not exist domains with boundary that lie in an open
great hemisphere of  $S^2$ and are wholly in one of the great
hemispheres of $S^3$ into which it is divided by the great sphere
$S^2$. In this case when $F$ is  a regular surface of class
$\mathcal{C}^2$, the saddle condition is equivalent to the
condition that the Gaussian curvature of $F^2$ does not exceed
one. We have the following result.
\begin{theorem}
\label{T2} \cite{Bor1}\cite{Bor2}\cite{Bor4}. Let $F$ be an
explicitly given compact saddle surface of smoothness class
$\mathcal{C}^1$ in the spherical space $S^3$. Then $F$ is a
totally geodesic great sphere.
\end{theorem}

 This theorem is a generalization of a
theorem of Bernshtein to a spherical space. For regular space we
obtain the following corollary.
\begin{theorem}\label{T3}\cite{Bor1}\cite{Bor2}
Let $F$ be an explicitly given compact surface that is regular of
class $\mathcal{C}^2$ in the spherical space $S^3$. If the
Gaussian curvature $K$ of $F$ satisfies $K\leqslant 1$ then $F$ is
a totally geodesic great sphere.
\end{theorem}
This theorem was stated in \cite{Bor2}. Really theorems
\ref{T2},\ref{T3} had been proved in \cite{Bor4} but were
formulated there for a centrally symmetric surfaces. The final
version was in \cite{Bor1}.

It seems to us that the following conjecture must hold under a
restriction on the Gaussian curvature of the surface. Suppose that
$F$ is an embedded compact surface, regular of class
$\mathcal{C}^2$, in the spherical space $S^3$. If the Gaussian
curvature $K$ of  $F$ satisfies $0< K\leqslant 1$, then $F$ is a
totally geodesic great sphere.

A.D. Aleksandrov \cite{Alex} had proved that an analytical surface
in Euclidean space $E^3$ homeomorphic to a sphere is a standard
sphere if principal curvatures satisfy the inequality
\begin{equation}\label{eq1}
(k_1+c)(k_2+c)\leqslant 0
\end{equation}
 This result had been generalized for analytic surfaces in
 spherical space $S^3$ and Lobachevsky space $L^3$ \cite{Bor4}:
 \begin{enumerate}
 \item [a)] in $S^3$ with additional hypothesis of positive
 Gaussian curvature;
 \item[b)] in $L^3$ under additional assumptions that principal
 curvatures $k_1, k_2$ satisfy $|k_1|$,\\$ |k_2|>c_0>1$.
 \end{enumerate}
 But in Lobachevsky space the result is true under  weaker analytic
 restriction.
 \begin{theorem}\label{T4}
  Let $F$ be a $\mathcal{C}^3$ regular surface homeomorphic to the
  sphere in the Lobachevsky space $L^3$. If $|k_1|, |k_2|>c_0>1$
  and principle curvatures $k_1$ and $k_2$ satisfy \eqref{eq1},
 then the surface is an umbilical round sphere in
  $L^3$.
 \end{theorem}
  Analogical result it is true for surfaces in the de Sitter space
  $S_1^3$.
  \begin{theorem}\label{T5}
  Let $F$ be a $\mathcal{C}^3$  regular compact spacelike surface
  in the de Sitter space $S^3_1$. If $|k_1|, |k_2|<1$ and
  principal curvatures satisfy \eqref{eq1}, then the
  surface is an umbilical round sphere in $S_1^3$.
  \end{theorem}

 Let $S_1^3$ be a simply-connected pseudo-Riemannian space of curvature
$1$   and signature $ (+, +,-)$. It can be isometrically embedded
in the pseudo-Euclidean space $E^4_{1}$  of signature
 $(+, +, +, -)$ as   the hypersurface given
 by the equation $x^2_1 +x_2^2 + x^2_3 - x^2_0 = 1$.  Together with $E^4_{1}$ we consider the superimposed Euclidean space $E^4$
  with unit sphere $S^3$ given by the equation  $x^2_1 +x_2^2 + x^2_3 + x^2_0 = 1$.
   We specify a mapping of $S^3_1$  into $S^3$. To the point $P$ of $S^3_1$
    with position vector $r$ we assign the point $\tilde P$ with position
    vector $\tilde{r} = {r}/\sqrt{1 + 2x^2_0}$.
     Under the mapping, to a surface $F \subset S^3_1$ corresponds a
    surface $\tilde{F}\subset S^3$. Let $b_{ij}$ and $\tilde{b}_{ij}$
     be the coefficients of the second quadratic
     forms of $F$ and $\tilde F$, and  $n=(n_1, n_2, n_3, n_0)$ be a normal vector field on $F$.
     \begin{lemma}\label{L1}\cite{Bor4}
     $\tilde{b}_{ij}={b_{ij}}/{\sqrt{1+2x_0^2}\sqrt{1+2n_0^2}}.$
     \end{lemma}
{\bf     Proof of theorem \ref{T1}.}
 From the condition $K\geqslant 1$ it follows that $F$ is a compact
 spacelike surface in the de Sitter space $S_1^3$. Locally a
 spacelike surface is explicitly given over totally geodesic
 great sphere $S_0^2\subset S_1^3$ and the orthogonal projection
 $p\colon F\to S_0^2$ in $S_1^3$ is covering. Indeed, $p$ is a
 local diffeomorphism. The compactness of $F$ and the simply
 connectedness of $S_0^2$ imply that $p$ is a global
 diffeomorphism $F$ on $S_0^2$ and the surface $F$ is globally
 explicitly given over $S_0^2$.

 We map a surface $F$ in $S_1^3$
into a surface $\tilde{F}$ in $S^3$. If $F$ has a definite metric
and Gaussian curvature $K\geqslant 1$, then $\tilde{F}$ has
Gaussian curvature  not greater than $1$. This follows immediately
from Lemma \ref{L1}, Gauss's formula  and the fact that $\langle
n, n\rangle = -1$ for normals to $F$. In a pseudo-Euclidean space,
the analogous  correspondence between surfaces and their
curvatures was used by Sokolov \cite{Sok}.

The surface $\tilde{F}$ satisfies the conditions of theorem
\ref{T3}. It follows that $\tilde{F}$ is a totally geodesic great
sphere. By lemma \ref{L1} the ranks of the second quadratic forms
of $\tilde{F}$ and $F$ coincide and we obtain that the surface $F$
is a totally geodesic surface in $S_1^3$.

{\bf Proof of Theorem \ref{T4} and \ref{T5}}. The normal $n(u_1,
u_2)$ to $F$ is chosen so that the principal curvature satisfy
\eqref{eq1}. In  a neighborhood of  an arbitrary nonumbilical
point $P$ we choose coordinate curves consisting of the  lines of
curvature, and an arbitrary orthogonal net in the case of
umbilical point. At $P$ the coefficients of the first quadratic
form are $e=g=1, f=0$. Let $F_1$ be the surface with radius vector
$\rho=(r-cn)/{\sqrt{|c^2-1|}}$.

In both cases the surface $F_1$ lies in $S_1^3$. Moreover
$$
\rho_{u_1}=\dfrac{(1+ck_1)}{\sqrt{|c^2-1|}}r_{u_1},\quad
\rho_{u_2}=\dfrac{(1+ck_2)}{\sqrt{|c^2-1|}}r_{u_2}.
$$
The unit normal $n_1=\dfrac{cr-n}{\sqrt{|c^2-1|}}$. From the
conditions on the principal curvatures of $F$ in
theorems~\ref{T4},~\ref{T5} it follows that
$$
\langle \rho_{u_1}, \rho_{u_1}\rangle>0,\quad \langle \rho_{u_2},
\rho_{u_2}\rangle>0
$$
 and $F_1$ is a spacelike
surface in $S_1^3$. The coefficients of the second quadratic form
of the surface $F_1$ are
$$
L_1=\dfrac{(1+ck_1)(k_1+c)}{\sqrt{c^2-1}},\quad
N_1=\dfrac{(1+ck_2)(k_2+c)}{\sqrt{c^2-1}}.
$$
The Gaussian curvature of $F_1$ at the point $P_1$ is equal to
$$
K=1-\dfrac{(k_1+c)(k_2+c)|c^2-1|}{(1+k_1c)^2(1+k_2c)^2}\geqslant1
$$
The same is true in umbilical points too. The surface $F_1$
satisfies the conditions of theorem \ref{T1}. It follows that the
surface $F_1$ is a totally geodesic great sphere in $S_1^3$ and
$F$ is an umbilical surface in $L^3$ or $S^3_1$.

\end{document}